\documentclass{amsart}
\usepackage{amsmath,amssymb,latexsym}
\usepackage[T1]{fontenc}
\usepackage{enumerate}

\def\FC{\mathrm{FC}}

\def\Ind#1#2{#1\setbox0=\hbox{$#1x$}\kern\wd0\hbox to 0pt{\hss$#1\mid$\hss}
\lower.9\ht0\hbox to 0pt{\hss$#1\smile$\hss}\kern\wd0}

\def\Notind#1#2{#1\setbox0=\hbox{$#1x$}\kern\wd0\hbox to 0pt{\mathchardef
\nn="3236\hss$#1\nn$\kern1.4\wd0\hss}\hbox to 0pt{\hss$#1\mid$\hss}\lower.9\ht0
\hbox to 0pt{\hss$#1\smile$\hss}\kern\wd0}

\newtheorem{theorem}{Theorem}[section]
\newtheorem{prop}[theorem]{Proposition}
\newtheorem{fact}[theorem]{Fact}
\newtheorem{lemma}[theorem]{Lemma}
\newtheorem{cor}[theorem]{Corollary}
\newtheorem*{claim}{Claim}

\theoremstyle{definition}
\newtheorem{defn}[theorem]{Definition}
\newtheorem{remark}[theorem]{Remark}

\def\pf{\par\noindent{\em Proof. }}
\def\pfclaim{\par\noindent{\em Proof of Claim. }}

\title{Centralizers in pseudo-finite groups}
\author{Nadja Hempel and Daniel Palac\'in}
\thanks{The second author was partially supported by the project MTM2014-59178-P}
\address{Department of Mathematics, University of California Los Angeles, Los Angeles, CA 90095-1555, USA}
\email{nadja@math.ucla.edu}
\address{Abteilung f\"ur Mathematische Logik, Mathematisches Institut,
	Albert-Ludwig-Universit\"at Freiburg, Ernst-Zermelo-Stra\ss e 1, D-79104
	Freiburg, Germany}
\email{palacin@math.uni-freiburg.de}

\begin{document}

\begin{abstract}
The role of finite centralizers of involutions in pseudo-finite groups is analyzed. It is shown that a pseudo-finite group admitting a definable involutory automorphism fixing only finitely many elements is finite-by-abelian-by-finite. As a consequence, we give a model-theoretic proof of a result for periodic groups due to Hartley and Meixner. 
Furthermore, it is shown that any pseudo-finite group has an infinite abelian subgroup.
\end{abstract}

\maketitle

\section{Introduction}

Certain structural properties of finite groups are inherited from the size of the centralizer of an involutory automorphism. This phenomenon was first noticed by Burnside who pointed out that a finite group admitting a fixed-point-free involutory automorphism is abelian. Another example is a relevant theorem of Brauer and Fowler who proved that there are only finitely many simple non-abelian finite groups containing an involution with a centralizer of a given size. Fong \cite{Fong} showed that if a finite group contains an involution whose centralizer has size $n$, then the group has a normal solvable subgroup whose index is bounded by a function on $n$. In fact, he proved a similar statement for elements of prime order assuming the Classification of Finite Simple Groups. Using Fong's Theorem, Hartley and Meixner \cite{HarMei} showed that any finite group $G$ admitting an involutory automorphism fixing only $n$ many elements has a normal nilpotent subgroup of class two whose index in $G$ only depends on $n$. In a series of papers the latter result has been extended first to arbitrary prime orders in \cite{HarMei2} and later to more general settings and even to arbitrary orders. We refer the interested reader to \cite{Shum} for clear exposition. 

The techniques to prove these style of results are typically from finite group theory. On the other hand, in Section $2$ of this paper we shall use techniques of infinite groups to show the aforementioned result of Hartley and Meixner, with some additional remarks on definable properties (see Theorem \ref{ThmMain1}). More precisely, we first prove this result for {\em pseudo-finite groups}, some model-theoretic limit groups, and then deduce the finite (and even the periodic) group version by a standard argument, yielding an alternative proof. Furthermore, using this approach we obtain in Corollary \ref{CorBelSes} that if $G$ is a group of exponent at most $m$ and admits an involutory automorphism $\alpha$ fixing only $n$ many elements, then $[G,\alpha]$ has index at most $n$ and its derived subgroup has a finite size which only depends on $n$ and $m$. This refines a result of Belyaev and Sesekin \cite{BelSes} on periodic groups, see also \cite{Bel}. 

Formally, a pseudo-finite group is an infinite group which satisfies every first-order sentence in the pure language of groups that is true of all finite groups.  Typical examples of pseudo-finite groups are torsion-free divisible abelian groups, infinite extra special groups of exponent $p>2$ and rank $n$, and linear groups over a pseudo-finite field. On the other hand, any free group forms an example of an infinite group which is not pseudo-finite. Equivalently, a pseudo-finite group $G$ is a group that satisfies the same elementary properties in the pure language of groups as a non-principal ultraproduct of an infinite family $\{G_i\}_{i\in \mathbb N}$ of finite groups. Therefore, by \L os's Theorem, which roughly speaking asserts that a first-order sentence in the pure language of groups is true of $G$ if and only if it is true of almost all finite groups $G_i$, one can see pseudo-finite groups as a model-theoretic limit of finite groups. For an overview on basic background around ultraproducts and known (model-theoretic) results on pseudo-finite groups we refer the reader to \cite{MacP}.

In Sections $3$ and $4$ we use finite group theory to deduce structural properties of pseudo-finite groups. More precisely, as in the locally finite case, using the Feit-Thompson Theorem we prove in Theorem \ref{ThmMain2} that any pseudo-finite group contains an infinite abelian subgroup. 
On the other hand, adapting some arguments of Shalev from \cite{Shalev}, which involve the Classification of Finite Simple Groups, we show in Theorem \ref{ThmMain3} that any pseudo-finite group  in which the centralizer of any element is finite or has finite index, the FC-center is a finite index definable subgroup. Under some strong model-theoretic assumptions, some results towards the existence of infinite abelian subgroups were already obtain in \cite{EJMR}. Moreover, Wagner in \cite{Wagner} recently has shown the existence a wide finite-by-abelian subgroup for pseudo-finite groups satisfying a weak chain condition on centralizers. 

Throughout the paper certain familiarity with basic notions on model theory and pseudo-finite constructions is assumed. In particular, we shall use without mentioning that any definable section of a pseudo-finite group is also pseudo-finite. A proof of this can be found in \cite[Lemma 2.9]{OuldPoint}.

\section{Involutions in pseudo-periodic groups}

By a {\em pseudo-periodic} group we mean an infinite group which is elementary equivalent to a non-principal ultraproduct of distinct periodic groups. In particular, any pseudo-finite group is obviously pseudo-periodic. 

In this section we analyze structural properties of pseudo-periodic groups containing involutions.
We start by pointing out an easy lemma on involutions in pseudo-periodic groups, which is proven easily using a well-known trick due to Brauer.

\begin{lemma}\label{LemInv}
	In any pseudo-periodic group, for any two involutions $i$ and $j$ either there is an involution centralizing both of them or there is some non-trivial element $z$  such that $j^z=i$ and $z^2=ji$.
\end{lemma}
\pf In any periodic group, for any two involutions, say $i$ and $j$, either $ij$ has even order $2m$ and hence $(ij)^m$ is a central involution of $\langle i,j\rangle$, or otherwise $ij$ has odd order $2m+1$ and so $j^{(ij)^m}=i$. Thus, either there is an involution centralizing $i$ and $j$, or there is some non-trivial element $z$ such that $j^z=i$ and $z^2=ji$. As the latter is first-order expressible in the pure language of groups, the same holds in any pseudo-periodic group. \qed

Recall that the {\em FC-center} of a group $G$ is the set of elements $x$ of $G$ whose conjugacy class $x^G$ is finite, or equivalently that the index $[G:C_G(x)]$ is finite. It is clearly a characteristic subgroup of $G$ and we denote it by $\FC(G)$. Using a result of Neumann \cite{Neu} stating that an infinite group cannot be covered by finitely many cosets of subgroups of infinite index, Belyaev  \cite[Lemma 2.1]{Bel} show the following:

\begin{fact}\label{FactFc}
	Suppose that $X$ is a finite subset of an arbitrary group $G$ and $X\cap X^g$ is non-empty for any $g$ in $G$. Then $X\cap \FC(G)$ is non-empty as well.
\end{fact}

An inspection of the proof of \cite[Lemma 2.2]{Bel} yields the following result.  

\begin{lemma}\label{LemFCInv}
	Let $G$ be a pseudo-periodic group containing an involution $i$ with a finite centralizer. Then either $\FC(G)$ contains an involution centralizing $i$ or $G=C_G(i)\cdot\FC(G)$.
\end{lemma}
\pf For any element $g$ of $G$, let $X_g$ be the set of involutions that commute with $i$ or $i^g$, and let $Y_g$ denote the finite set
$$
X_g \cup g^{-1} C_G(i) \cup C_G(i)g^{C_G(i)}. 
$$
We aim to apply Fact \ref{FactFc} to the set $Y_g$. To do so, we prove that for any element $x$ of $G$ we have that $Y_g^x \cap Y_g$ is non-empty. Suppose that there exists some element $x$ in $G$ such that $X_g^x \cap X_g$ is empty. In particular, this yields that there is no involution commuting with $i^x$ and $i$ nor with $i^g$ and $i^x$. Hence by the previous lemma there are two elements $y$ and $z$ of $G$ such that $i^{gy} = i^x$ with $y^2=i^gi^x$ and $i^{xz}=i$ with $z^2=i^xi$.
Thus, we trivially have that $y^{-1} i^g = i^x y^{-1}$ and $z^{-1} i^x = i z^{-1}$,
and so
$$
(y^{-1}z^{-1})^{i^x} = i^x y^{-1} z^{-1} i^x = y^{-1} i^g i z^{-1} = y^{-1} i^g i^x i^x i z^{-1} = y^{-1} y^2 z^2 z^{-1} = yz.
$$
Moreover, as $i^x = z i z^{-1}$, we have that $(y^{-1}z^{-1})^{i^x} = ( z^{-1} y^{-1} )^{iz^{-1}}$ 
and hence $iz^{-1}$ inverts $(yz)^{-1}$. On the other hand, as $i^{gyz} = i$ we get that $(yz)^{-1}$ belongs to $C_G(i)g$ and so there is an element $a$ in $C_G(i)$ such that $(yz)^{-1}=ag$. In addition, notice that $xz$ and so $xzi$ also commute with $i$. Thus, there exists some $b$ such that $iz^{-1} = bx$. Therefore, putting all together we have that $(ag)^{bx} = (ag)^{-1}$, so
$$
(C_G(i)g^{G_G(i)})^x \cap g^{-1}C_G(i) \neq \emptyset,
$$
and hence $Y_g^x \cap Y_g$ is never empty. 

Consequently, by Fact \ref{FactFc} the intersection $Y_g \cap \FC(G)$ is non-empty for any element $g$ of $G$. Thus, either there is an element $g$ in $G$ for which $X_g\cap \FC(G)$ is non-empty, yielding the existence of an involution in $\FC(G)\cap C_G(i)$, or for every $g$ in $G$ the intersection of $\FC(G)$ with $g^{-1}C_G(i)\cup C_G(i)g^{C_G(i)}$ is non-empty implying that every $g$ belongs to $C_G(i)\cdot \FC(G)$.  \qed

Recall that a group is said to be {\em definably simple} if it has no definable normal subgroup in the pure language of groups. As an immediate consequence of the fact above we deduce the following:

\begin{prop}
	There is no infinite pseudo-periodic definably simple group containing an involution $i$ with a finite centralizer.
\end{prop}
\pf If such a group exists, it is not abelian and so it has a trivial center. Furthermore, Lemma \ref{LemFCInv} yields the existence of a non-trivial element $x$ such that $C_G(x)$ has finite index in $G$. Hence, as this centralizer is a proper subgroup and $x$ has only finitely many conjugates, the definable normal proper subgroup $C_G(x^G)$ of $G$ has also finite index and thus it is non-trivial, a contradiction. \qed

As an application to the previous results, we can give now an easy proof of a version of the Brauer-Fowler lemma for definably simple periodic groups, shown by Shunkov \cite{Shunkov} in the simple periodic case.

\begin{cor}
	There are only finitely many periodic definably simple groups containing an involution with a centralizer of a given finite size, and all these groups are finite.
\end{cor}
\pf Observe first that any definably simple periodic group containing an involution with a finite centralizer must be a finite simple group by the previous result (or its proof using \cite[Lemma 2.2]{Bel}). Now, suppose towards a contradiction that there are infinitely many finite simple groups with a centralizer of an involution of size $n$, and let $G$ an ultraproduct of all of them with respect to some non-principal ultrafilter. This is clearly an infinite definably simple pseudo-finite group. However, such an infinite group contains an involution with a finite centralizer, contradicting the previous result.  \qed

In fact, Shunkov \cite{Shunkov} showed that any periodic group containing an involution with a finite centralizer is locally finite (and moreover solvable-by-finite). Using this result, Hartley and Meixner \cite{HarMei} proved the following:

\begin{fact}\label{FactHarMei}
	There exists a function $f(n)$ defined on the natural numbers such that if $G$ is a periodic group containing an involution $i$ with $|C_G(i)|\le n$, then $G$ contains a nilpotent subgroup of class two and index at most $f(n)$.
\end{fact}

Using that such a periodic group is locally finite, the proof reduces then to the finite case by a classical argument of Kegel and Wehrfritz \cite[page 54]{KW}. Consequently,  their proof uses techniques from finite group theory as well as Fong's result \cite{Fong}. In fact, Hartley \cite[Theorem 1.1]{Har1} pointed out that an easy argument yields the following:

\begin{fact}\label{FactHar}
	There is a function $f(n)$ defined on the natural numbers such that if $G$ is a periodic group admitting an involutory automorphism $\alpha$ with $|C_G(\alpha)|\le n$, then $G$ contains a normal subgroup $H$ of $G$ of index bounded by $f(n)$ and whose derived subgroup is contained in $C_G(\alpha)$.
\end{fact}

On the other hand, Belyaev and Sesekin \cite{BelSes} using methods typical for infinite groups obtained the following result, see also \cite{Bel} for some additional information on the FC-center of the group. 

\begin{fact}\label{FactBelSes}
	Let $G$ be a periodic group admitting an involutory automorphism $\alpha$ such that $C_G(\alpha)$ is finite. Then $[G,\alpha]'$ and $G/[G,\alpha]$ are finite.
\end{fact}

In particular, as an immediate corollary one obtains that the centralizer of $[G,\alpha]'$ is a finite index nilpotent subgroup of class two. However, the proof does not provide information on the index of this nilpotent subgroup as in fact it does not give a bound on the size of the derived subgroup depending only on $n$. 

Next, we give an alternative proof of the aforementioned results of Hartley and Meixner, obtaining  a gently improvement regarding definability properties in Fact \ref{FactHarMei} and \ref{FactHar}, as well as bounds in Fact \ref{FactBelSes} for groups of bounded exponent.

First we prove some basic lemmata on pseudo-finite groups.

\begin{lemma}\label{LemCenFin}
	Let $G$ be a pseudo-finite group containing an element $x$ with a finite centralizer, and let $N$ be a definable normal subgroup of $G$. Then 
	$$
	|C_{G/N}(xN)| \le |C_G(x)|.
	$$
	Moreover, if $x$ is an involution and the equality holds, then $N$ is contained in and centralizes the set $xx^G$. In particular, it is a subgroup of $Z([x,G])$ and is inverted by $x$.
\end{lemma}
\pf As $G$ is pseudo-finite, there is an infinite family $\{G_i\}_{i\in\mathbb N}$ of finite groups such that $G$ is elementary equivalent to an ultraproduct of this family, for some non-principal ultrafilter on $\mathbb N$. Thus without loss of generality we may assume that $G$ is this ultraproduct. Let $N_i$ denote the set given by the formula defining $N$ in $G_i$, which in almost all such groups is a normal subgroup, and let $x_i$ be the element of $G_i$ corresponding to a representative of $x$. Observe that for almost all finite groups $G_i$ we have that $C_{G_i}(x_i)$ has the same size as $C_G(x)$, which is finite by assumption. Additionally, since for finite groups it holds that
$$
\frac{ |C_{G_i}(x_i/N_i)| }{|N_i|} = |C_{G_i/N_i}(x_iN_i) | \le | C_{G_i}(x_i) |,
$$
see for instance \cite[Lemma 1.1]{HarMei}, we deduce that $N_i$ has index at most $|C_G(x)|$ in $C_{G_i}(x_i/N_i)$, and so does $N$ in $C_G(x/N)$ as the latter is an elementary property. It then follows that $C_{G/N}(xN)$ is of size at most $|C_G(x)|$, as desired. 

For the second part of the statement, notice that if the equality holds, then so does in almost all finite groups $G_i$ and consequently we have that 
$$
\frac{ |C_{G_i}(x_i/N_i)| }{ |C_{G_i}(x_i)| } = |N_i|.
$$
Thus, considering the surjective map $y\mapsto [x_i,y]$ from $C_{G_i}(x_i/N_i)$ to $N_i$, we easily get that each element of these groups $N_i$ is of the form $x_ix_i^y=[x_i,y]$ for some $y$ in $G_i$. Hence $N_i$ is contained in the set $x_ix_i^{G_i}$ and whence the same holds replacing $N_i$, $G_i$ and $x_i$ by $N$, $G$ and $x$. Consequently, the involution $x$ inverts every element of $N$. Moreover, since $N$ is normal, any conjugate of $x$ inverts $N$ and hence the latter centralizes the set $xx^G$, as desired.  \qed

The lemma above is well-known for finite groups or even for arbitrary groups when $N$ is assumed to be finite, see \cite[Lemma 1.3]{Bel}. Similarly, the results below follow easily from the finite setting. 

\begin{lemma}\label{LemInvComm}
	Let $G$ be a pseudo-finite group admitting a definable automorphism $\alpha$ 
	with $|C_G(\alpha)|\le n$ for some natural number $n$. Then the subgroup $[G,\alpha]$ is definable of index at most $n$.
\end{lemma}
\pf 
As $G$ is pseudo-finite, there is a family $\{G_m\}_{m\in\mathbb N}$ of finite groups whose ultraproduct $\mathbb G$ with respect to some non-principal ultrafilter is elementary equivalent to the infinite group $G$. In particular, almost all groups $G_m$ contain a definable automorphism $\alpha_m$ with a centralizer of size $n$, which is induced by the same first-order formula that defines $\alpha$. Moreover, we denote the definable automorphism in $\mathbb G$ still as $\alpha$. Set $X_m$ to be the set of elements of the form $x^{-1}x^{\alpha_m}$ with $x$ in $G_m$, and note that in almost all finite groups $G_m$ we trivially have that
$$
\frac{ |G_m| }{ |X_m |} = |C_{G_m}(\alpha_m)| \le n.
$$
Let $X$ be the set of elements of the form $x^{-1}x^{\alpha}$ with $x$ in $\mathbb G$ which generates the subgroup $[\mathbb G,\alpha]$ of $\mathbb G$. Now, suppose the group $[\mathbb G,\alpha]$ has index at least $n+1$ in $\mathbb G$. Then we can find $n+1$ translates of $X$ which are pairwise disjoint. Consequently, this has to be true in almost all finite groups for the corresponding $X_m$ which contradicts the above inequality. Thus the group $[\mathbb G,\alpha]$ has index at most $n$ in $\mathbb G$. Moreover, this subgroup is the union of all $k$-products of elements in $X \cup X^{-1}\cup \{1\}$ for any natural number $k$. Thus, by saturation and compactness this subgroup is indeed equal to the $k$-products of elements in $X \cup X^{-1}\cup\{1\}$ for some fixed natural number $k$.  This remains true for $G$ and so the subgroup $[G,\alpha]$ is definable in $G$ and has as well index at most $n$.
\qed

A standard argument using ultraproducts yields the following:

\begin{cor}\label{CorDefCom}
	There is a function $f(n)$ defined on the natural numbers such that if $G$ is a locally finite group admitting an involutory automorphism $\alpha$ with $|C_G(\alpha)|\le n$ for some natural number $n$, then the subgroup $[G,\alpha]$ has index at most $n$ and is obtained as the product of at most $f(n)$ many elements of the form $x^{-1}x^\alpha$ with $x$ in $G$.
\end{cor}
\pf Let $G$ be locally finite and assume, as we may after replacing $G$ by $G\rtimes\langle\alpha\rangle$, that $\alpha$ is an element of $G$. It is immediate that $[G,\alpha]$ has index at most $n$: Assume there exist  elements $x_0,\ldots,x_n$ of different cosets of $[G,\alpha]$ in $G$. Let $H$ be the finite subgroup of $G$ generated by the elements $\alpha,x_0,\ldots,x_n$. We obtain that 
$$
[ H : [H,\alpha] ] \le C_H(\alpha) \le n.
$$ 
which clearly  contradicts the choice of $H$.

For the second assertion, we first find such a function for finite groups. To do so, suppose towards a contradiction that there is no such a function. Thus, there is a natural number $n$ such that for any natural number $k$ we can find a finite group $G_k$, admitting an involutory automorphism $\alpha_k$ fixing only $n$ many elements, in a way that in $[G_k,\alpha_k]$ there is an element which cannot be written as the product of $k$ many elements of the form $x^{-1}x^{\alpha_k}$. Now, consider an ultraproduct $K$ of all finite groups $G_k\rtimes \langle \alpha_k\rangle$ with respect to some non-principal ultrafilter on $\mathbb N$. Notice that $K$ is a pseudo-finite group such that the equivalence class  $\alpha$ of the sequence $((1,\alpha_k))_{k\in\mathbb N}$ has a centralizer of size $n$. Thus, the previous lemma yields that $[K,\alpha]$ is definable and so, by saturation of $K$ and compactness, it is indeed defined as the product of at most $l$ many elements of the form $x^{-1}x^\alpha=[x,\alpha]$ for some natural number $l$. As saying that the product of $l+1$ of such commutator elements is indeed the product of $l$ many is first-order expressible, \L os' Theorem yields that in almost all subgroups $G_k\rtimes \langle \alpha_k\rangle$ the  commutator subgroup $[G_k\rtimes \langle \alpha_k\rangle,(1,\alpha_k)]$ is also defined as the product of at most $l$ many commutator elements. Therefore, as an easy computation gives us that
$$
[G_k\rtimes \langle \alpha_k\rangle,(1,\alpha_k)] = [G_k,\alpha_k] \times \{1\},
$$
we obtain an isomorphism between the latter and $[G_k,\alpha_k]$, which contradicts the choice of $G_k$. 

To finish, notice that the same function obtained in the finite context yields the statement for locally finite groups. \qed

\begin{lemma}\label{LemFinNilp}
	Let $G$ be a pseudo-finite group and let $N$ be a definable central subgroup of $G$. If $G/N$ is finite-by-abelian then $G$ has a definable finite index subgroup $K$ that is nilpotent of class two. Moreover, if the derived subgroup of $G/N$ has size $n$, then the index of $K$ in $G$ only depends on $n$.
\end{lemma}
\pf Assume, as we may after replacing $G$ by the finite index subgroup $C_G(G'/N)$, that $G/N$ is in addition nilpotent of class two. Now, as $G$ is pseudo-finite, we can find an infinite family $\{G_i\}_{i\in\mathbb N}$ of finite groups in a way that $G$ is elementary equivalent to a non-principal ultraproduct of this family. Let $N_i$ denote the set given by the formula defining $N$ in $G_i$, which in almost all finite groups $G_i$ is a central subgroup. Moreover, as $(G/N)'=G'N/N$ is finite of size $n$ by assumption, almost all finite groups $G_i/N_i$ have also a derived subgroup of size $n$. 

Now, fix one of these finite groups $G_i$ and set $H_i=Z(G_i) \cap G_i' =C_{G_i'}(G_i)$. Note that 
$$
[G_i':H_i] \le [G_i':G_i'\cap N_i] = [G_i'N_i : N_i] = n.
$$
Thus, there are at most $n$ many elements $x_{i,1},\ldots,x_{i,n}$ of $G_i$ such that 
$$
H_i = C_{G_i'}(x_{i,1},\ldots,x_{i,n}).
$$
Now, set $K_i$ to be the subgroup $C_{G_i}(x_{i,1},\ldots,x_{i,n}/N_i)$. As $G_i'N_i/N_i$ has size at most $n$, so does the conjugacy class of $x_{i,j}N_i$ in $G_i/N_i$ and hence the subgroup $C_{G_i}(x_{i,j}/N_i)$ has index at most $n$ in $G_i$. Consequently, we get $[G_i:K_i]\le n^n$. Finally, set $F_i$ to be the subgroup of $G_i$ generated by $x_{i,1},\ldots,x_{i,n}$. Thus
$$
[[F_i,K_i],K_i] \le [G_i'\cap N_i, K_i] = \{1\}
$$ by construction and as $N_i$ is central. Therefore $[K_i',F_i]$ is trivial by the three subgroup lemma, yielding that $K_i'$ is contained in $C_{G_i'}(F_i)=H_i$. Since $H_i$ is central in $G_i$ we obtain that $K_i$ is nilpotent of class two. 

To conclude, observe that the index of $K_i$ in $G_i$ is at most $n^n$ and moreover, there is a first-order formula which only depends on $n$ that defines $K_i$. Hence, such a formula defines a subgroup $K$ of $G$ which has index at most $n^n$ and is also nilpotent of class two.  \qed

\begin{lemma}\label{LemNil2}
	Let $G$ be a pseudo-finite group admitting a definable involutory automorphism $\alpha$ with $|C_G(\alpha)|\le n$ for some natural number $n$ and let $N$ be a definable $\alpha$-invariant subgroup of finite index. If $N$ is nilpotent of class two, then there is a definable finite index subgroup $H$ of $N$ which is normal in $G$ and such that $H'$ is contained in $C_G(\alpha)$.
\end{lemma}
\begin{proof}
	Consider the group $[N,\alpha]$, a definable subgroup of $N$ of index at most $n$ by Lemma \ref{LemInvComm}. As $\alpha$ inverts any element of the form $x^{-1}x^\alpha$ with $x$ in $N$, for any two elements $u$ and $v$ of these form we have that
	$$
	[u,v]^\alpha = [u^\alpha,v^\alpha] = [u^{-1},v^{-1} ] = [v,u^{-1}]^{v^{-1}} = [ u , v ]^{u^{-1}v^{-1}} = [u,v]
	$$ 
	since $[N,\alpha]'$ is central in $N$. Moreover, as the latter is generated by the commutator elements of elements of the form $x^{-1}x^\alpha$, we deduce that $[N,\alpha]'$ is fixed pointwise by $\alpha$. Hence, to conclude the proof it is enough to take the intersection of the finitely many conjugates of $[N,\alpha]$ in $G$, which yields a definable normal subgroup with the desired properties. 
\end{proof}

Now we state and prove Hartley and Meixner's result for pseudo-finite groups.

\begin{prop}\label{PropMain}
	Let $G$ be a pseudo-finite group containing an involution $i$ with a finite centralizer. Then there is a definable normal subgroup $H$ of $G$ of finite index such that $H'$ is contained in $C_G(i)$.
\end{prop}
\pf 
First note that it is enough to show that $G$  is finite-by-abelian-by-finite: \\
In this case, there exists a finite index normal subgroup $A$ of $G$ which is finite-by-abelian. Suppose that the derived subgroup of $A$ has size $k$. Then, as the conjugacy class of any element of $A$ has size at most $k$, the FC-center of $G$ contains $A$ and so it has finite index in $G$ as well. Moreover, if $y_1,\ldots,y_m$ are representatives of the cosets of $A$ in $\FC(G)$, then every element of $\FC(G)$ has a conjugacy class of size at most $\max\{ k\cdot |y_i^G| \, : \, i=1,\ldots,m \}$
and so the $\FC$-center of $G$ is definable and finite-by-abelian by \cite[14.5.11]{Rob}. It follows that the definable normal subgroup $\FC(G)\cap C_{G}(\FC(G)')$ has finite index in $G$ and is nilpotent of class two. Thus by Lemma \ref{LemNil2} we obtain the result.

Now, let $G$ be a pseudo-finite group with an involution $i$ with a finite centralizer of size $n$ and suppose that any pseudo-finite group with an involution that has a centralizer of size strictly smaller than $n$ is finite-by-abelian-by-finite. 

Note that, using minimality with respect to the size of the centralizer of an involution and Lemma \ref{LemCenFin} (or \cite[Lemma 1.3]{Bel}) we may assume that any finite normal subgroup $H$ of $G$ is contained in and centralizes the set $ii^G$. In particular, any finite normal subgroup is abelian and inverted by $i$.

Let $T$ be the set of periodic elements of $\FC(G)$, which by a result of Neumann, see \cite[14.5.9]{Rob}, forms a characteristic subgroup of $G$. 
Suppose first that $T$ is finite. Then $\FC(G/T)$ has no elements of finite order and thus it is abelian as its derived subgroup is periodic again by \cite[14.5.9]{Rob}. Moreover, Lemma \ref{LemFCInv}  yields that $\FC(G/T)$ has finite index in $G/T$ and thus $G/T$ is abelian-by-finite and so $G$ is finite-by-abelian-by-finite. Thus, we may assume that $T$ is infinite. Then, the periodic subgroup $T$ is additionally locally normal-finite by Dicman's Lemma, see \cite[14.5.7]{Rob}. Consequently, by our assumption on all finite normal subgroups we obtain that $T$ is contained in and centralized by $ii^G$. Therefore, it is abelian and inverted by $i$. Now, let $N$ be a maximal normal subgroup of $G$ contained in $Z([i,G])$. Note that this is an infinite abelian subgroup and so it cannot contain $i$ as $C_G(i)$ is finite. Observe as well that $N$ is inverted by $i$. Moreover, as it centralizes the subgroup $[i,G]$, the subgroup $C_G(N)$ has finite index in $G$ by Lemma \ref{LemInvComm} (we use this fact again at the end of this proof). Hence, the latter equals to a finite intersection of centralizers of elements from $N$, whence it is definable and so is $Z(C_G(N))$. 

Set $\bar G$ to be $G$ modulo $Z(C_G(N))$. If this group is finite, then $G$ is abelian-by-finite and we can conclude. So we may assume that $\bar G$ is infinite. Furthermore, note that $Z(C_G(N))$ is an infinite abelian group and thus it does not contain $i$. Hence  $\bar G$ is an infinite pseudo-finite group containing an involution $\bar i$ with a centralizer of size $|C_{\bar G}(\bar i)| \le n$ by Lemma \ref{LemCenFin}. 

\begin{claim}  The group $\bar G$ is finite-by-abelian-by-finite.
\end{claim}
\pfclaim
If $|C_{\bar G}(\bar i)| <n $, we obtain that $\bar G$ is finite-by-abelian-by-finite by our initial assumption.
Thus we may assume that  $|C_{\bar G}(\bar i)| = n$. So Lemma \ref{LemCenFin} yields that the definable normal subgroup $Z(C_G(N))$ is contained in the set $ii^G\cap C_G(ii^G) \subseteq Z([i,G])$. As $N$ is contained in $Z(C_G(N))$, the choice of $N$ yields that $N=Z(C_G(N))$ and so $N$ is definable. As pointed out before, the elements $\bar T$ of finite order of $\FC(G/N)$ form a normal subgroup. If $\bar T$ is finite, arguing as before we obtain that $\bar G/\bar T$ is abelian-by-finite and hence $\bar G$ is finite-by-abelian-by-finite. If  $\bar T$ is infinite, using the same arguments as for $T$, we get that $\bar T$ is an infinite locally  finite-normal abelian subgroup. Moreover it is inverted by $i$ and contained in $Z([iN,G/N])$. Observe that $[iN,G/N]$ equals to $[i,G]$ modulo $N$ since $N$ is a normal subgroup of $[i,G]$, and so the preimage $P$ of  $\bar T$ in $G$ is contained in $[i,G]\cap C_G([i,G]/N)$. Furthermore, as $N$ is central in $[i,G]$, we deduce that $P$ is nilpotent of class two. 
As $\bar T$ is infinite and locally finite-normal we can find a subgroup $H$ of $P$ which is normal in $G$ and such that $H/N$ is finite but has size at least $n+1$. In particular, since $N$ is definable so is the subgroup $H$. As $H$ is a subgroup of $P$, we have that $H$ is trivially contained in $[i,G]\cap C_G([i,G]/N)$, is nilpotent of class two and $H/N$ is inverted by $i$.

Now, we show that $H$ modulo its central subgroup $C_{H}(ii^G)$ has exponent two. To do so, consider an element $h$ of $H$ and an element $g$ of $ii^{G}$. As $H$ is contained in $C_{G}([i,G]/N)$, there is some $x$ in $N$ such that $g^h=gx$. Moreover, as $H/N$ is inverted by $i$, there is some element $v$ of $N$ such that $h^{i}=h^{-1}v$. Note that $i$ inverts $x$ and $v$ as these belong to $N$ and also $g$ since it belongs to $ii^{G}$. Thus, since $v$ commutes with $h$ and $g$, we obtain that
$$
g^{hi} = (h^{i})^{-1} g^{i} h^{i} = v^{-1} h  g^{-1} h^{-1} v = h g^{-1} h^{-1}.
$$
On the other hand, as $x$ commutes with $g$, we have that
$$
g^{hi} = (gx)^{i} = g^{-1}x^{-1} = (xg)^{-1} = (gx)^{-1} = h^{-1} g^{-1} h
$$
and so $h^2$ commutes with $g$. As $g$ and $h$ were taken to be arbitrary we conclude. 

Thus, as $H/N$ is inverted by $i$ and $N$ is included in $C_{H}(ii^{G})$, we deduce that the group $H/C_{H}(ii^{G})$ is centralized by $i$ since all its elements are involutions. As $C_{H}(ii^{G})$ is central in $H$, the same holds for $H/Z(H)$. As $Z(H)$ is definable, the infinite group $G/Z(H)$ is pseudo-finite. Thus, applying Lemma \ref{LemCenFin} to $G/Z(H)$ we obtain that $|C_{G/Z(H)}(iZ(H))|\leq n$ and consequently $H/Z(H)$ has size at most $n$. Moreover if equality holds in the above equation, Lemma \ref{LemCenFin} yields that $Z(H)$ is a normal abelian subgroup of $G$ inverted by $i$ and contained in $Z([i,G])$. Hence $Z(H)$ is equal to $N$ by the choice of the latter, and so this yields that the size of $H/N$ in $\bar G$ is at most $n$, a contradiction. Hence $
|C_{G/Z(H)}(iZ(H))|<n.$ Whence, we deduce that $G/Z(H)$ is a pseudo-finite group with an involution that has a centralizer of size strictly smaller than $n$ and thus by our assumption it is finite-by-abelian-by-finite. Since $H$ is a finite extension of $N$ and $N$ is central in $H$, we get that  $\bar G = G/N$ is finite-by-abelian-by-finite as claimed. 
\qed$_{\mbox{Claim}}$

Using the claim, we can find a definable normal finite index subgroup $F$ of $G$ such that $F/Z(C_G(N))$ is finite-by-abelian. As mentioned before the claim, the subgroup $C_G(N)$ is definable and  has finite index in $G$. Thus taking the intersection of these two finite index subgroups, namely the subgroup $C_F(N)$, we obtain a finite index subgroup of $G$ which is finite-by-abelian modulo its center. Thus by Lemma \ref{LemFinNilp} $C_F(N)$ contains a definable subgroup $K$ of finite index that is nilpotent of class two. Consequently $K$ has finite index in $G$ and therefore has only finitely many conjugates in $G$. Considering the intersection of these conjugates, we obtain a definable finite index subgroup of $G$ which is nilpotent of class two and is additionally normal in $G$. Using Lemma \ref{LemNil2} we obtain the desired result. 
\qed

As a consequence we obtain the following result for periodic groups, from which standard arguments yield a gently improvement of Fact \ref{FactHarMei} and \ref{FactHar}, offering at the same time an alternative proof. 

\begin{theorem}\label{ThmMain1}
	There exists a function $f(n)$ defined on natural numbers such that if $G$ is a periodic group admitting an involutory automorphism $\alpha$ with $|C_G(\alpha)| \le n$, then there is a characteristic definable subgroup $H$ of $G$ such that $G/H$ and $H'$ have both size  bounded by $f(n)$.
\end{theorem}
\pf  
We first prove the statement for finite groups. Suppose, towards a contradiction, that such a function does not exist for finite groups. Thus, there exists a natural number $n$ such that for any natural number $k$ we can find a finite group $H_k$ admitting an involutory automorphism $\alpha_k$ that fixes only $n$ many elements of $H_k$ such  that each $H_k$ witnessing that $f(n)$ cannot be equal to $k$.  Now, for each natural number $k$, let $G_k$ be the semi-direct product $H_k\rtimes \langle \alpha_k\rangle$, a finite group containing an involution with a centralizer of size $n$, and let $G$ be an ultraproduct of all these finite groups with respect to some non-principal ultrafilter on $\mathbb N$. Notice that $G$ is an infinite pseudo-finite group containing an involution with centralizer of size $n$ and hence, Proposition \ref{PropMain} yields the existence of a definable normal subgroup $K$ of $G$ of finite index whose derived subgroup has size at most $n$. It then follows that $\FC(G)$ contains $K$ and so it has finite index in $G$. Furthermore, the FC-center is finite-by-abelian as explained at the beginning of the proof of Proposition \ref{PropMain}. 

Now, let 
$m= \max \{[G: \FC(G)], |\FC(G)'|\}.$
For each natural number $k$, consider the first-order formula $\phi_k(x)$ in the pure language of groups defining in an arbitrary group the set of elements $x$ with only $k$ many conjugates, or equivalently the set of elements $x$ whose centralizer has index at most $k$. 
By choice of $m$, we have that the formula $\phi_m(G)$ defines the characteristic subgroup $\FC(G)$ of $G$ whose derived group has size at most $m$.
Notice that having a derived group of size at most $m$ as well as having index at most $m$ is first-order expressible in the pure language of groups.
Thus the formula $\phi_{m}(x)$ induces a characteristic subgroup in almost all groups $G_k$ whose index and the size of its derived subgroup are bounded by $m$. Hence, it also induces a characteristic subgroup in $H_k$ with the same properties, contradicting the choice of the $H_k$.

For the general case, let $f(n)$ be the function given in the finite context and set $k$ to be $f(n)$. Now, suppose towards a contradiction that  there is some natural number $n$ such that for any natural number $m$ we can find a periodic group $G_m$ in a way that the infinite family $\{G_m\}_{m\in\mathbb N}$ witnesses a counterexample. Consider a non-principal ultraproduct $G$ of these periodic groups.

As each $G_m$ is locally finite by Shunkov's result, there is a directed system of finite groups $\{G_{m,l}\}_{l\in I}$ whose union equals $G_m$. It follows from the finite case that $k$ many translates of $\phi_k(G_{m,l})$ covers $G_{m,l}$ for every $l$ and $m$. Thus, we trivially have that 
$$
\frac{ |G_{m,l}| }{ |\phi_k(G_{m,l}) |} \le k.
$$
Now, suppose the subgroup $\FC(G)$ has index at least $k+1$ in $G$. Then we can find $k+1$ translates of $\phi_k(G)$ which are pairwise disjoint and therefore, this has to be true in almost all infinite groups $G_m$ by \L os' Theorem. Thus in each $G_m$ we can then find elements $x_0,\ldots,x_k$ witnessing that $\phi_k(G_m)$ has $k+1$ pairwise disjoint translates. Now, as $G$ is the union of the directed system of   finite groups $\{G_{m,l}\}_{l\in I}$, there is $l\in I$ which contains the elements $x_0,\ldots,x_k$. This contradicts the above inequality and thus $[G:\FC(G)]$ is bounded by $k$, which only depends on $n$. 

Now, as $\FC(G)$ is given as the union of the subsets $\phi_m(G)$ and the group $G$ is $\aleph_0$-saturated, an easy compactness argument, together with a similar argument as before, yields the existence of a natural number $m_0$ such that 
$\phi_{m_0}(x)$ induces a characteristic subgroup in almost all periodic subgroups $G_m$ whose index and the size of its derived subgroup are both bounded by $m_0$, a contradiction. \qed

\begin{remark} By Theorem \ref{ThmMain1}, one easily gets the statements of Fact \ref{FactHarMei} and \ref{FactHar}. Namely, by taking $C_H(H')$ one obtains a characteristic  definable nilpotent subgroup of class two whose index only depends on $n$, yielding the definable version of Fact \ref{FactHarMei}. Furthermore, similar arguments as in Lemma \ref{LemNil2}, together with Corollary \ref{CorDefCom}, gives the definable version of Fact \ref{FactHar}. 
\end{remark}

Observe that our methods do not yield a first-order formula only depending on $n$ that defines uniformly these characteristic subgroups. As in the proof of Belyaev and Sesekin of Fact \ref{FactBelSes}, the main reason is that we cannot control the size of the derived subgroup of the FC-center by a function depending on $n$. Nevertheless, we obtain a uniform bound for groups of bounded exponent.

\begin{cor}\label{CorBelSes}
	There is a function $f(n)$ defined on natural numbers  such that for any group $G$  of exponent at most $m$ admitting an involutory automorphism $\alpha$ with $|C_G(\alpha)| \le n$, the subgroup $[G,\alpha]$ has index at most $n$ in $G$ and its derived subgroup has size at most $f(n,m)$.
\end{cor}
\pf By Shunkov's Theorem, after Corollary \ref{CorDefCom}, it remains to prove that the size of such a derived subgroup can be bounded by a function only depending on $n$ and $m$. As usual, suppose towards a contradiction that this is not the case. Thus, we can find two natural numbers $n$ and $m$ and a family $\{H_k\}_{k\in\mathbb N}$ of groups of exponent $m$ admitting an involutory automorphism $\alpha_k$ fixing only $n$ many elements such that $[H_k,\alpha_k]'$ has size at least $k$. Now, consider the ultraproduct $G$ of all groups $H_k\rtimes\langle \alpha_k\rangle$ with respect to some non-principal ultratrafilter on $\mathbb N$. Thus, the infinite group $G$ is again a group of exponent at most $m$ and the equivalence class $\alpha$ of the sequence $(\alpha_k)_{k\in\mathbb N}$ is an involution with a centralizer of size $n$. Hence, by Theorem \ref{ThmMain1} we obtain a finite-by-abelian normal subgroup of finite index and so the FC-center of $G$ has also finite index in $G$. 

Now, we prove that $[G,\alpha]$ is finite-by-abelian. To do so, let $L$ denote the intersection of all conjugates of $[G,\alpha]\cap \FC(G)$, a finite index normal subgroup of $G$. Moreover, since $G$ is periodic, its FC-center is locally finite-normal, see \cite[14.5.7]{Rob}, and so $L$ is the union of finite normal subgroups of $G$. By Lemma \ref{LemCenFin} there exists a finite subgroup $N$ of $L$, normal in $G$ such that the size of $C_{G/N}(\alpha N)$ is minimal among all finite subgroups of $L$ which are normal in $G$. Hence, again by Lemma \ref{LemCenFin} we obtain that $L/N$ is contained in $Z([G/N,\alpha N])$. As $L$ is a normal subgroup of $[G,\alpha]$ by construction, we obtain that $[G,\alpha]$ modulo $N$ is central-by-finite and so finite-by-abelian by a result of Schur, see \cite[10.1.4]{Rob}. It then follows that $[G,\alpha]$ is also finite-by-abelian since $N$ is finite, as claimed.

By Corollary \ref{CorDefCom}, there is a first-order formula $\psi(x,y)$ in the pure language of groups such that $\psi(x,\alpha)$ defines $[G,\alpha]$ and also that $\psi(x,(1,\alpha_k))$ defines 
$$
[H_k\rtimes\langle\alpha_k\rangle,(1,\alpha_k)] = [H_k,\alpha_k]\times \{1\}.
$$
Moreover, as expressing that the derived subgroup of the group defined by $\psi(x,\alpha)$ has size at most $|[G,\alpha]'|$ is first-order expressible, the same is true of $\psi(x,(1,\alpha_k))$ in almost all the semi-direct products $H_k\rtimes\langle\alpha_k\rangle$ and so $[H_k,\alpha_k]'$ has size $|[G,\alpha]'|$. However, this yields a contradiction, finishing the proof. \qed

To finish this section, we point out that a standard elementary argument allows us to extend the previous results to the pseudo-periodic context, which extends Proposition \ref{PropMain}. As the proofs are similar to the one of Theorem \ref{ThmMain1} we omit to give details.

\begin{theorem}
	There is a function $f(n)$ defined on natural numbers such that if $G$ is a pseudo-periodic group admitting a definable involutory automorphism $\alpha$ with $|C_G(\alpha)| \le n$, then there is a characteristic definable subgroup $H$ of $G$ such that $G/H$ and $H'$ have both size  bounded by $f(n)$.
\end{theorem}
\pf Firstly, notice that Theorem \ref{ThmMain1} yields the existence of a function $f_0$ such that any such pseudo-periodic group is covered by at most $f_0(n)$ many translates of the subset of elements with at most $f_0(n)$ many conjugates. Hence, if the statement does not hold, one can find an infinite family of counter-examples and consider a non-principal ultraproduct $G$. Notice that $G$ is also  covered by $f_0(n)$ many translates of its corresponding subset of elements with at most $f_0(n)$ many conjugates. Hence, a similar argument as in the proof of Theorem \ref{ThmMain1} yields that the FC-center of $G$ is given by a first-order formula and it is finite-by-abelian, which contradicts the choice of our family. \qed 

Consequently, following the proof of Lemma \ref{LemNil2} and using Corollary \ref{CorDefCom} we obtain the following:

\begin{cor}
	There is a function $f(n)$ defined on natural numbers such that if $G$ is a pseudo-periodic group admitting a definable involutory automorphism $\alpha$ with $|C_G(\alpha)| \le n$, then there is a definable normal subgroup $H$ of $G$ whose index is bounded by $f(n)$ and its derived subgroup is contained in $C_G(\alpha)$.
\end{cor}

\section{The existence of an infinite abelian group}

In this section we prove the following:

\begin{theorem}\label{ThmMain2}
	Any pseudo-finite group has an infinite abelian subgroup.
\end{theorem}
\pf
We first show that any pseudo-finite group has a non-trivial element with an infinite centralizer. To do so, suppose towards a contradiction that all non-trivial elements of a pseudo-finite group $G$ have finite centralizer.
Thus the group $G$ cannot contain an involution by Lemma \ref{LemFCInv} as otherwise the FC-center would be non-trivial.
On the other hand, by the Feit-Thompson Theorem, any finite group $H$ without an involution is solvable. Thus, in any such finite group there exists a non-trivial element $h$ such that $\langle h^H \rangle$ is abelian and so, the conjugacy class $h^H$ is contained in $C_H(h)$. Hence, as the latter is first-order expressible, the same holds in the pseudo-finite group $G$. Namely, we can find a non-trivial element $g$ in $G$ whose conjugacy class is included in its centralizer. As by assumption, the centralizer of $g$ is finite, so is $g^G$ and thus $G$, a blatant contradiction. Therefore, we deduce that any pseudo-finite group contains a non-trivial element with an infinite centralizer.

Now, as an element of infinite order yields an infinite abelian subgroup, we may assume that $G$ is periodic. Let $x_0$ be an element of $G$ with infinite centralizer and consider the pseudo-finite group $G_0$ defined as $C_G(x_0)$ modulo $\langle x_0\rangle$. Again, we can find an element $x_1$ in $C_G(x_0)$ whose class $\bar x_1$ in $G_0$ is non-trivial and has an infinite centralizer. As $x_0$ and $x_1$ commute and both have finite order, the group $\langle x_0,x_1 \rangle$ is finite. Thus, considering the pseudo-finite group $G_1$ defined as $C_G(x_0,x_1)$ modulo $\langle x_0,x_1 \rangle$, we find an element $x_2$ in $C_G(x_0,x_1)$ whose class $\bar x_2$ in $G_1$ is non-trivial and has again an infinite centralizer. Proceeding in this way, we obtain an infinite abelian subgroup generated by the set $\{x_i\}_{i\in\mathbb N}$. This finishes the proof. \qed

As an immediate consequence of the theorem above we obtain the following:

\begin{cor}\label{CorAbelian}
	There is only a finite number of finite groups in which every abelian subgroup has size at most $n$.
\end{cor}

\section{Restricted Centralizers}

In this section we explore groups with a certain condition on centralizers, weakening the notion of FC-group, introduced by Shalev \cite{Shalev}.

\begin{defn}
	A group has {\em restricted centralizers} if the centralizer of any element is either finite or of  finite index.
\end{defn}

A subgroup of a group with restricted centralizers has also restricted centralizers. Moreover, Shalev showed in \cite[Lemma 2.1]{Shalev} that this condition is preserved under taking quotients by finite normal subgroups. Observe  that for pseudo-finite groups, having restricted centralizers is preserved for any definable factor by Lemma \ref{LemCenFin}.

Model-theoretically, natural examples of groups with restricted centralizers are groups of $\mathrm{SU}$-rank $1$, such as infinite extra-special groups of exponent $p>2$, which is also pseudo-finite. On the other hand, Tarski monsters also satisfy such property, witnessing that a priori this condition may not be strong enough to ensure a well-behaved group theoretic structure. Nevertheless, Shalev showed \cite[Theorem 1.1]{Shalev} that any  profinite group with restricted centralizers is finite-by-abelian-by-finite. His proof makes use of the Classification of Finite Simple Groups as well as Zelmanov's solution of the Burnside problem for periodic profinite groups.
We obtain the corresponding result for pseudo-finite groups, which require the following result obtained combining \cite{HarMei2} and \cite[Theorem 2]{Khu}:

\begin{fact}
	There are two functions $f(n,p)$ and $h(p)$ defined on the natural numbers such that if $G$ is a finite group admitting an automorphism of prime order $p$ fixing only $n$ many elements, there is a nilpotent subgroup of index bounded by $f(n,p)$ and whose nilpotency class is at most $h(p)$.
\end{fact}

\begin{theorem}\label{ThmMain3}
	The FC-center of any pseudo-finite group with restricted centralizers is definable and has finite index.
\end{theorem}
\pf Let $G$ be a pseudo-finite group having restricted centralizers and assume, as we may, that it is not an FC-group. By \cite[Lemma 2.2]{Shalev}, the group $G$ modulo $\FC(G)$ is periodic and so there exists some element $x$ of $G\setminus \FC(G)$ such that $x^p$ belongs to $\FC(G)$ for some prime number $p$. Thus, the definable subgroup $C_G(x^p)$ has finite index in $G$ and, since by assumption $C_G(x)$ is finite, say of size $n$, it additionally admits a definable automorphism of order $p$ which fixes only $n$ many elements.

Now, as the group $C_G(x^p)$ is pseudo-finite, it is elementary equivalent to the ultraproduct of some family $\{H_m\}_{m\in \mathbb N}$ of finite groups. Notice that almost all these finite groups admit an automorphism of order $p$ which fixes $n$ many elements. By the previous fact, in each of these $H_m$ we can find a nilpotent subgroup of index bounded by $f(n,p)$, for some function, and whose nilpotency class only depends on $p$. As such a nilpotent subgroup has at most $f(n,p)$ many conjugates, the intersection of all these yields a normal nilpotent subgroup of the same nilpotency class whose index is bounded in terms of $n$ and $p$. Namely, its index is at most $f(n,p)!$. It then follows that the Fitting subgroup $F(H_m)$ of each $H_m$ has also index at most $f(n,p)!$ and its nilpotency class is bounded by a function on $n$ and $p$. In particular, for almost all natural numbers $m$, the nilpotency class of $F(H_m)$ and the index of the latter subgroup in $H_m$ are the same. On the other hand, there is a first-order formula defining in a uniform way all these Fitting subgroups and hence, this formula yields also the existence of a normal nilpotent subgroup $N$ of $C_G(x^p)$ of finite index. Furthermore, after replacing $N$ by the intersection of its finitely many conjugates we may assume that it is normal in $G$.

To conclude, as $N$ has finite index in $G$, it is enough to prove that $\FC(G)$ extends $N$. To do so, it suffices to reproduce the proof of \cite[Proposition 2.5]{Shalev} where it is shown that any normal nilpotent subgroup of a group with restricted centralizers is contained in the FC-center of $G$. This finishes the proof. \qed

In case that the pseudo-finite group is assumed to be $\aleph_0$-saturated, then an easy compactness argument yields that its FC-center is finite-by-abelian. For finite groups the result above yields the following:

\begin{cor}\label{CorResCen}
	There is a function $f(n)$ defined on the natural numbers such that if $G$ is a finite group such that for any element $x$ the size of $C_G(x)$ or $G/C_G(x)$ is at most $n$,
	then there is a characteristic subgroup $H$ of $G$ such that the size of $G/H$ and $H'$ are bounded by $f(n)$.
\end{cor}
\pf If the assertion does not hold, then for some natural number $n$ we can find an infinite set $\{G_k\}_{k\in\mathbb N}$ of finite groups satisfying the condition such that $G_k$ has no characteristic subgroup $H_k$ with $|G/H_k|$ and $|H_k'|$ bounded by $k$. Thus, any ultraproduct $G$ of these finite groups with respect to a non-principal ultrafilter on $\mathbb N$ has restricted centralizers and so its FC-center is definable and has finite index. Hence, again by compactness and $\aleph_0$-saturation of $G$ we get that $\FC(G)$ is indeed finite-by-abelian and so $G$ is finite-by-abelian-by-finite. However, the formula defining $\FC(G)$ induces in almost all finite groups $G_k$ a characteristic subgroup $H_k$ such that $|G_k/H_k| = |G/\FC(G)|$ and also that $|H_k'| \le |\FC(G)'|$, a contradiction.\qed

Note that a pseudo-finite group of $\mathrm{SU}$-rank $1$ has uniform restricted centralizers and so satisfies the weak chain condition on centralizers. Thus, we deduce that these groups are finite-by-abelian-by-finite (cf. \cite[Corollary 4.12]{Wagner} and \cite[Theorem 2.11]{EJMR}). In fact, this can be deduced from  Theorem \ref{ThmMain2} since any infinite abelian subgroup is contained in the definable finite-by-abelian subgroup.

\end{document}